\documentclass{amsart}

\usepackage{amssymb}
\usepackage{amsmath}
\usepackage{amsthm}
\usepackage{amsbsy}

\theoremstyle{plain}
\newtheorem{theorem}{Theorem}[section]
\newtheorem{lemma}[theorem]{Lemma}
\newtheorem{corollary}[theorem]{Corollary}

\theoremstyle{definition}
\newtheorem{definition}[theorem]{Definition}

\theoremstyle{remark}

\newcommand{\C}{\mathcal{C}}
\newcommand{\N}{\mathbb{N}}

\begin{document} 

\title{The derived series and virtual Betti numbers} 

\author{Siddhartha Gadgil}

\address{Stat-Math Unit,\\
	Indian Statistical Institute,\\
	Bangalore, India}
\email{gadgil@isibang.ac.in}

\subjclass{Primary 57M05, 57M10; Secondary 57M60}

\begin{abstract}
The virtual Betti number conjecture states that any hyperbolic
three-manifold has a finite cover with positive first Betti number.
We show that this would follow if it were known that the derived
series of the fundamental group $G$ of a hyperbolic three-manifold
satisfies a certain stability property. The stability property is the
statement that if all the quotients $G^i/G^{i+1}$ of the derived
series $G^i$ of $G$ are finite, then the derived series stabilises.

The proof involves basic facts regarding finite group actions on
homology spheres.

\end{abstract}

\date{\today}

\maketitle

An important conjecture in hyperbolic geometry is the so called
virtual Betti number conjecture which states that any hyperbolic
$3$-manifold $M$ has a finite cover with positive first Betti
number. A natural way to try to construct such a cover is to consider
the derived series $G^i$ of $G=\pi_1(M)$. There are \emph{a priori}
three possibilities for this: the series may stabilise, some quotient
$G^i/G^{i+1}$ may be infinite or we may have finite derived quotients
without the series stabilising. In the second case we obtain a finite
cover with positive Betti number.

Our goal in this note is to show that to prove the virtual Betti
number conjecture for a class (closed under passing to finite covers)
of hyperbolic $3$-manifolds, it suffices to rule out the third case
for manifolds in the class. We first recall some basic definitions
and elementary facts.

Suppose $G$ is a (finitely presented) group. The \emph{derived series}
of $G$ is defined inductively as follows: let $G^0=G$ and, for
$i\in\N$, let $G^i=[G^{i-1},G^{i-1}]$. We say that the derived series
stabilises if $G^{k+1}=G^k$ for some $k$. In this case, for $i\geq k$
we have $G^i=G^k$. The quotient $G^i/G^{i+1}$ is just the
abelianisation of $G^i$ and hence the derived series stabilises if and
only if, for some $i$, $G^i$ is perfect.

Suppose $G=\pi_1(M)$ for a hyperbolic $3$-manifold $M$. Then the
virtual Betti number conjecture is equivalent to the statement that
there is a finite-index subgroup $H$ of $G$ whose abelianisation is
infinite. Suppose $G^i/G^{i+1}$ is infinite for some $i\geq 0$, then
if we let $H=G^i$ for \emph{the smallest} such $i$, $H$ has finite
index in $G$ and has infinite abelianisation as required.

To state our result, we need a definition.

\begin{definition}
Let $\C$ be a class of hyperbolic $3$-manifolds that is closed under
passing to finite covers. Then we say that $\C$ satisfies the
\emph{stability property} if for each $M\in\C$ the following holds for
$G=\pi_1(M)$: if each quotient $G^i/G^{i+1}$, $i\in\N$, of the derived
series of $G$ is finite, then the derived series stabilises.
\end{definition}

We could take the class $\C$ to be, for instance, all arithmetic
manifolds or simply all hyperbolic $3$-manifolds. Our main result is
the following.

\begin{theorem}
Suppose the class $\C$ of hyperbolic $3$-manifolds satisfies the
stability property. Then any manifold $M\in\C$ has a finite cover with
positive first Betti number.
\end{theorem}

To prove the result, we need the following lemma.

\begin{lemma}
Let $M$ be a hyperbolic $3$-manifold and let $N$ be a finite Galois
cover of $M$ with $H$ the group of deck transformations. Suppose
$H_1(M)=H_1(N)=0$, then $H$ is the binary icosahedral group (or $M=N$).
\end{lemma}
\begin{proof}
By Poincar\'e duality, $M$ and $N$ are homology spheres. As $H_1(M)=0$
and $H$ is a quotient of $\pi_1(M)$, it follows that $H$ is
perfect. Thus $H$ is a finite perfect group acting freely on the
homology $3$-sphere $N$ (in particular $H$ has periodic
cohomology). It is well known (see, for instance,~\cite{Mi}
or~\cite{Wo}) that $H$ must be the binary icosahedral group (or the
trivial group).
\end{proof}

In particular, we have the following corollary as the binary
icosahedral group has order~$120$.

\begin{corollary}
Let $M$ be a hyperbolic $3$-manifold and let $N$ be a finite Galois
cover of $M$ with order more than $120$. Then $N$ is not a homology
sphere.
\end{corollary}

Suppose now that $\C$ is a class of manifolds satisfying the
hypothesis of the theorem and let $M$ be a manifold in $\C$. Consider
the covers of $M$ corresponding to the derived subgroups of
$\pi_1(M)$. By the stability property $M$ either has a finite cover
with positive first Betti number or $M$ has a cover that is a homology
sphere. Thus we may assume that $H_1(M)=0$.

By residual finiteness, we can find a finite Galois cover $M'$ of $M$
that has order greater than $120$. Consider the derived series of $M'$
and the corresponding covers. By the stability condition we either get
a finite cover $N$ with positive first Betti number or we have a
finite cover $N$ with $H_1(N)=0$. It suffices to rule out the latter
case.

Assume that we are in the latter case and that $N$ is a finite cover
of $M'$ with $H_1(N)=0$ corresponding to a derived subgroup of
$\pi_1(M')$. As the derived subgroups are characteristic subgroups of
$\pi_1(M')$ which is in turn a normal subgroup of $\pi_1(N)$,
$\pi_1(N)$ is a normal subgroup of $\pi_1(M)$. Thus $N$ is a Galois
cover of $M$. As the degree of the covering map is greater than $120$,
we get a contradiction to the above corollary. This completes the
proof of the theorem.

\qed

\bibliographystyle{amsplain}

\end{document}